\documentclass[12pt]{amsart}
\input{epsf}
\usepackage{amscd,amssymb}

\newtheorem{thm}{Theorem}[section]
\newtheorem{cor}[thm]{Corollary}
\newtheorem{lemma}[thm]{Lemma}
 \newtheorem{prop}[thm]{Proposition}

\theoremstyle{definition}

\newtheorem{ex}{Example}

\theoremstyle{remark}

\newcommand{\R}{\mathbb R} 
\newcommand{\affine}{\mathbb E}
\newcommand{\E}{\mathbb E}
\newcommand{\B}{\mathbb B}

\newcommand{\Isom}{\operatorname{Isom}}

\newcommand{\oto}{\operatorname{O}(2,1)}
\newcommand{\soto}{\operatorname{SO}(2,1)}

\newcommand{\vx}{{\mathsf x}}

\newcommand{\vv}{{\mathsf v}}

\newcommand{\tx}{p_o}
\newcommand{\ty}{p_1}
\newcommand{\tz}{p_2}
\newcommand{\tu}{q_o}
\newcommand{\tv}{q_1}
\newcommand{\tw}{q_2}

\newcommand{\xp}[1]{\vx^+(#1)}
\newcommand{\xm}[1]{\vx^-(#1)}
\newcommand{\xo}[1]{\vx^o(#1)}
\newcommand{\xpm}[1]{\vx^\pm(#1)}

\newcommand{\bg}[1]{{\mathcal B}_{#1}}
\newcommand{\ggp}[1]{{\mathcal G}^+_{#1}}
\newcommand{\ggm}[1]{{\mathcal G}^-_{#1}}

\renewcommand{\gg}{\gamma}
\newcommand{\gG}{\Gamma}
\newcommand{\ga}{\alpha}
\newcommand{\gl}{\lambda}
\newcommand{\gt}{\tau}

\newcommand{\cT}{\mathcal{T}}
\newcommand{\cL}{\mathcal{L}}
\newcommand{\cS}{\mathcal{S}}
\newcommand{\cF}{\mathcal{F}}

\newcommand{\origin}{\mathsf 0}

\newcommand{\cqfd}{q.e.d.}

\begin{document}

\title{Closed timelike curves in flat Lorentz spacetimes}
\author{Virginie Charette}
    \address{Department of Mathematics\\ McMaster University\\
    Hamilton, Ontario, Canada}
    \email{charette@math.mcmaster.ca}
\author{Todd A.\ Drumm}
\address{Department of Mathematics and Statistics \\
    Swarthmore College\\ Swarthmore, PA\ \ 19081, USA}
    \email{tad@swarthmore.edu}
\author{Dieter Brill}
    \address{Department of Physics \\
    University of Maryland\\ College Park, MD\ \ 20742, USA}
    \email{brill@physics.umd.edu}
\date{\today}

\begin{abstract}
We consider the region of closed timelike curves ({\em CTC}'s) in
three-dimensional flat Lorentz spacetimes. The interest in this
global geometrical feature goes beyond the purely mathematical.
Such spacetimes are lower-dimensional toy models of sourceless
Einstein gravity or cosmology. In three dimensions all such
spacetimes are known: they are quotients of Minkowski space by a
suitable group of Poincar\'e isometries. The presence of CTC's
would indicate the possibility of ``time machines", a region of
spacetime where an object can travel along in time and revisit the
same event. Such spacetimes also provide a testbed for the
chronology protection conjecture, which suggests that quantum back
reaction would eliminate CTC's. In particular, our interest in
this note will be to find the set free of CTC's for $\E/<\gamma>$,
where $\E$ is modeled on Minkowski space and $\gamma$ is a
Poincar\'e transformation. We describe the set free of CTC's where
$\gamma$ is hyperbolic, parabolic, and elliptic. \end{abstract}
\maketitle

Let $\E$ denote three-dimensional {\em Minkowski space}.  This is
an affine space with translations in $\R^{2,1}$, the vector space
equipped with the standard indefinite bilinear form
$\B(\cdot,\cdot)$ of signature $(2,1)$. (Since $\E$ is flat and
geodesically complete, the reader may identify $\E$ with its set
of translations without any major difficulties arising.) We look
at flat Lorentz manifolds, which are quotients of an open subset
$X$ of $\E$ by a group $\Gamma$ of affine Lorentzian isometries
that acts properly discontinuously on that subset. Such manifolds
$X/\Gamma$ inherit a local causal structure from $\E$.

Let $X/\Gamma$ be a flat Lorentz manifold. A {\em timelike} vector
is a vector $\vv\in\R^{2,1}$ such that $\B(\vv,\vv)<0$.  A {\em
timelike curve} in $X/\Gamma$ is a $C^{1}$ path
$c:[0,1]\rightarrow X/\Gamma$ whose tangent vectors are all
timelike; we say that $c$ is {\em closed} if $c(0)=c(1)$. The
purpose of this note is to lay the groundwork for understanding
regions free of closed timelike curves in a Lorentz spacetime.

The {\em CTC region} of a spacetime is the set of all points which
lie on some closed timelike curve and the {\em CTC-free region} is
the complement of this space. Suppose $\Gamma$ acts properly
discontinuously on some subset $X\subset\E$. Denote by $[p]$ the
image of $p$ in $X/\Gamma$ under projection. We wish to determine
the set of all $p$ such that $[p]$ lies in the CTC region of $X
/\Gamma$.

We first note the following basic lemma.
\begin{lemma}\label{lem:timelikecurves} Let $\Gamma$ and
$X\subset\affine$ be as above.  Let $p\in X$ be a point such that
$\gg (p)-p$ is a timelike vector, for some $\gg\in\Gamma$, and the
line segment starting at $p$ and ending at $\gg (p)$ lies entirely
in $X$. Then $[p]\in X/\Gamma$ lies on a smooth closed timelike
curve. \end{lemma} This lemma is well known, but we provide a
proof for completeness. \proof Let $c:\R\rightarrow\E$ be a
continuous and piecewise linear path through each point
$\gg^{n}(p)$, defined by the following: \begin{equation*}
c(t)=\gg^{[t]}(p)+(t-[t])\left( \gg^{[t]+1}(p)
-\gg^{[t]}(p)\right), \end{equation*} where $[t]$ denotes the
largest integer no greater than $t$. The path $c$ is of the class
$C^{\infty}$ except at the integers. We describe a new path near
$t=0$ and then apply the same procedure at all integer values of
$t$. For some small $\epsilon >0$, define \begin{equation*}
\tilde{c}(t)=d_{0,\epsilon}(t)\left( p+ t(p-\gg^{-1}(p)\right)+
u_{0,\epsilon}(t) \left( p+ t(\gg(p)-p)\right) , \end{equation*}
where $d_{0,\epsilon}(t)$ and $u_{0,\epsilon}(t)$ are $C^{\infty}$
functions such that \begin{equation*} \begin{array}{lcr}
d_{0,\epsilon}(t) =
    \left\{ \begin{array}{l} 1 \mbox{ for } t\leq 0
            \\ 0 \mbox{ for } t\geq \epsilon
            \end{array} \right. &
\mbox{ and } & u_{0,\epsilon}(t)
    = \left\{ \begin{array}{l} 0 \mbox{ for } t\leq 0
            \\ 1 \mbox{ for } t\geq \epsilon . \end{array}
    \right.
    \end{array}
\end{equation*}
The path $\tilde{c}$ is $C^{\infty}$ and agrees with the path $c$
for $t\leq 0$ and $t \geq \epsilon$. The path $\tilde{c}$ can be
extended to a $C^{\infty}$ path through each point $\gg^{n}(p)$.
We can assume that $\gg^{n-1}(p)-\gg^{n}(n)$ is a future pointing
timelike vector. The sum of two future pointing timelike vectors
is another future pointing timelike vector. Therefore, $\tilde{c}$
is a smooth timelike curve which goes through each point
$\gg^{n}(p)$. \cqfd

\begin{cor}
Suppose $\Gamma$ acts properly discontinuously on $\E$.  Then
$[p]$ is in the CTC region of $\E/\Gamma$ if and only if
$\gamma(p)-p$ is timelike, for some $\gamma\in\Gamma$. \end{cor}
The vector $\gg(p)-p$ is called the {\em displacement vector} for
$p$.

In this note we will mainly consider $\Gamma =\langle \gg
\rangle$, where $\gamma$ is an affine Lorentzian isometry, that
is, an element of the {\em Poincar\'e group}.  Note that the
cyclic group $\Gamma$ acts properly discontinuously on all of $\E$
if and only if $\gamma$ admits no fixed points.  We will examine
both the fixed point case  -- in other words, when $\gg$ is an
element of the {\em Lorentz group} -- and the fixed point-free
case.

\section{Isometries}
Let $\Isom(\E)$ denote the group of affine isometries of $\E$,
that is the Poincar\'e group.  Choosing an origin $\origin$, we
may write the action of an isometry $\gamma\in\Isom(\E)$ as:
\begin{equation}\label{eq:basic}  \gg(p)=\origin + g(p-\origin)
+\vv, \end{equation} where $g\in \oto$ is called its {\em linear
part}  and $\vv\in\R^{2,1}$ is  called its {\em translational
part}. (Typically, we shall omit $\origin$, which can be taken to
be the point whose coordinates are all zero.)

Denote the Lorentzian inner product on the vector space of
translations $\R^{2,1}$ by $\B \left(  \cdot, \cdot \right)$:
\begin{equation*} \B \left( \begin{bmatrix}x\\y\\z\end{bmatrix},
\begin{bmatrix}u\\v\\w\end{bmatrix}\right)=xu+yv-zw.
\end{equation*}
It is invariant under the action due to any element of $\oto$.

The union of the sets of timelike vectors and non--zero lightlike
vectors divides into two connected components, one of which is
called {\em future directed} and one called {\em past directed}.
It is common to choose the connected component where the third
coordinate is positive to be the future direction.  Such a choice
is called a {\em time orientation} on $\E$. We say that a set
$X\subset\E$ is {\em future (resp. past) complete} if given a
point $p\in X$ every future (resp. past) directed ray starting at
$p$ remains in $X$.

We will restrict our examination to the identity component
$G\subset\Isom(\E)$, consisting of orientation and
time-orientation preserving isometries.  This is a subgroup of
$\soto$ which has index $2$.

The conjugacy classes of elements of $G$ are  identified by their
trace. Further, a nonidentity element $g\in G$ is called
\begin{itemize}
    \item {\em hyperbolic} if $tr(g)>3$,
    \item {\em parabolic} if $tr(g)=3$ and
    \item {\em elliptic} if $tr(g)<3$.
\end{itemize}

We say that an affine transformation is {\em hyperbolic}, {\em
parabolic}, or {\em elliptic} if its linear part is hyperbolic,
parabolic, or elliptic, respectively. Conjugacy classes of
hyperbolic and parabolic affine transformations are determined by
the trace of the linear part and what we will call the {\em signed
Lorentzian length}.  This invariant, due to
Margulis~\cite{Margulis1,Margulis2}, measures the Lorentzian
length of the closed geodesic determined by the isometry; more on
this in the following section.

For a given transformation $\gg$, $\E$ splits into three regions,
according to the causal character of the displacement vector
$\gamma(p)-p$, for $p\in\E$: \begin{itemize} \item the {\em
timelike region} for the action of $\gg$ on $\E$:
\begin{equation*} \cT(\gg)  = \{ p\in\E | \B\left(
\gg(p)-p,\gg(p)-p \right) <0\}; \end{equation*}

\item the {\em lightlike region} for the action of $\gg$: \begin{equation*}
\cL(\gg)  = \left\{ p\in\E | \B\left( \gg(p)-p,\gg(p)-p \right) =
0 \right\} ; \end{equation*}

\item the {\em spacelike region } for $\gg$:
\begin{equation*}
\cS(\gg)  =  \{ p\in\E | \B\left( \gg(p)-p,\gg(p)-p )\right. > 0
\} . \end{equation*} \end{itemize}

Suppose $\gG \subset\Isom(\E)$ acts properly discontinuously on
the maximal set $X\subset\E$; we denote by $\cT(\gG)$ the set of
points $p\in X$ such that $[p]$ is in the CTC region of $X/\gG$.
If $\gG$ acts freely and properly discontinuously on $\E$, then:
\begin{equation*} \cT(\gG) = \bigcup_{\gg\in\gG} \cT(\gg) .
\end{equation*} The region whose quotient is free of closed
timelike curves is denoted by $\cF(\gG)$. If $\gG$ acts freely and
properly discontinuously on $\E$, then: \begin{equation*} \cF(\gG)
= \bigcap_{\gg\in\gG} \left(\cS(\gg) \cup \cL(\gg) \right).
\end{equation*}

\subsection{Conjugation}
An isometry may be conjugated by an element of $GL(3,\R)$, in
order to facilitate explicit calculations.  For instance, given a
hyperbolic $g\in G$, there exists $h\in GL(3,\R)$ such that
$hgh^{-1}$ is a diagonal matrix.

Conjugation of an isometry by a linear map corresponds to a change
of basis in $\E$. Conjugation by a pure translation corresponds to
changing the origin.

Conjugation of an isometry by elements of the Poincar\'e group is
more restrictive but has real physical meaning. Conjugation by an
element of $\soto$ corresponds to changing to an observer in a
different inertial reference frame. Conjugation by a translation
corresponds to changing to an observer at a different space-time
location. The Principle of Special Relativity requires the
invariance of physical laws under a change of location and
inertial reference.

The relevant invariants in our Lorentzian spacetimes, such as the
Lorentzian inner product and signed Lorentzian length, are
invariant after conjugation by an element of the Poincar\'e group.
In particular, even though the form of \eqref{eq:basic} depends on
the choice of origin, the signed Lorentzian length does not depend
on the choice of origin.

\section{Hyperbolic transformations}
Let  $\gg\in\Isom(\E)$ be a hyperbolic transformation with linear
part $g$.  Then $g$ admits three eigenvectors, $\vx^-(g),
\vx^o(g), \vx^+(g)$ with eigenvalues $\gl(g) < 1 < \gl(g)^{-1}$,
respectively. The eigenvectors $\vx^{\pm}(g)$ are null and can be
chosen so that they are future directed and $\B
(\vx^{-}(g),\vx^{+}(g) ) = -1$.  Then we may choose $\vx^o(g)$ to
be the unique 1-eigenvector satisfying: \begin{itemize} \item it
is unit-spacelike, i.e. $\B(\vx^o(g),\vx^o(g))=1$; \item
$\{\vx^-(g), \vx^+(g),\vx^o(g)\}$ is a right--handed basis
    for $\R^{2,1}$.
\end{itemize}

The {\em signed Lorentzian length} of $\gamma$ is:
\begin{equation*} \ga (\gg) = \B(\gg (p)-p,\vx^o(g)),
\end{equation*} where $p\in\affine$ is arbitrary (i.e.
$\alpha(\gamma)$ does not depend on the choice of $p$). Every
hyperbolic transformation $\gg$ gives rise to a unique invariant
line $C_{\gg}\subset\E$, which is parallel to $\vx^o(g)$. If $q\in
C_{\gg}$:
\begin{equation*}
\gg(q)=q+\ga(\gg)\vx^o(g);
\end{equation*}
thus $\ga (\gg)$ measures the signed Lorentzian length of
$C_\gg/\langle \gg \rangle$. The curve $C_\gg/\langle \gg \rangle$
is called the {\em unique closed geodesic} in $\affine/\langle \gg
\rangle$, that is $C_\gg/\langle \gg \rangle$ is the image of a
simple closed curve $c[0,1]\rightarrow\E$ such that $c'(t)$ is
constant on $[0,1]$, where $c'(0)$ is thought of as a right hand
limit and $c'(1)$ is thought of as a left hand limit.

As trace identifies conjugacy classes of elements in $\soto$, the
trace of the linear part and the value of $\ga$ identify conjugacy
classes of hyperbolic isometries in $\Isom(\E)$. Conjugate $\gg$
by the change of basis matrix $[ \vx^-(g) ~ \vx^+(g) ~ \vx^o(g) ]$
(an element of $GL(3,\R)$, not $\soto$), so that its linear part
is diagonal. Next, conjugate by an appropriate translation, so
that the origin lies on $C_\gg$.

Write the point $p\in\E$ in terms of the new basis and origin:
\begin{equation*} p= p_- \vx^-(g) + p_+ \vx^+(g) + p_o \vx^o(g).
\end{equation*} We may write the conjugate transformation as
follows: \begin{equation*}
\begin{bmatrix}p_- \\ p_+ \\ p_o \end{bmatrix} \longmapsto
          \begin{bmatrix} \gl(g) & 0 & 0 \\
                            0 & \gl(g)^{-1} & 0 \\
                            0 & 0 & 1 \end{bmatrix}
          \begin{bmatrix} p_- \\ p_+ \\ p_o \end{bmatrix}
          + \begin{bmatrix} 0 \\ 0 \\ \ga(\gg) \end{bmatrix}
          =   \begin{bmatrix}\gl(g) p_-\\ \gl(g)^{-1}p_+\\ p_o+\ga(\gg)
          \end{bmatrix}.
\end{equation*}
The causal character of the displacement vector is determined by
the following calculation: \begin{equation*} \B(\gg
(p)-p,\gg(p)-p)= 2 \B(\vx^+(g), \vx^-(g))(\gl(g)-1)(\gl(g)^{-1}
-1)p_{-}p_{+} +\ga(\gg)^2. \end{equation*} Thus points in
$\cT(\gg)$, i.e. points projected onto CTC's, satisfy the
following inequality:
\begin{equation}\label{eq:hyperbolicspacelike} p_{-} p_{+} <
\frac{-\ga(\gg)^2 }{2 (1- \gl(g))(\gl(g)^{-1} -1)}. \end{equation}
Observe that the right hand side is negative.

\subsection{Hyperbolic transformations with fixed points}
Suppose that $\gamma$ admits a fixed point; equivalently,
$C_\gg$ is pointwise fixed and $\ga(\gg)=0$.
Set:
\begin{align*} \bg{\gg}&=\{ p_-\xm{\gg}+ p_+\xp{\gg} +
    p_o\xo{\gg} ~|~ p_-p_+<0\} .
\end{align*}
Then $p$ solves \eqref{eq:hyperbolicspacelike} if and only if
$p\in\bg{\gg}$. This region divides into two connected components,
bounded by $\cL(\gg)$, which in turn is composed of the so-called
{\em unstable}/{\em stable} planes $E^{\pm}(\gg)$: these are the
planes containing $C_{\gg}$ and parallel to $\vx^{\pm}(g)$.

\begin{figure}
    \caption{A cross section for hyperbolic transformations}
    \label{fig:crosssection}
\begin{picture}(220,220)
    \put(110,110){\vector(1,1){80}} \put(110,110){\vector(-1,-1){80}}
    \put(110,110){\vector(-1,1){80}} \put(110,110){\vector(1,-1){80}}

    \put(195,195){\makebox(0,0){$E^{+}(\gg)$}}
    \put(15,195){\makebox(0,0){$E^{-}(\gg)$}}

    \put(110,165){\makebox(0,0){$\ggp{\gg}$}}
    \put(110,55){\makebox(0,0){$\ggm{\gg}$}}

    \put(55,110){\makebox(0,0){$\bg{\gg}$}}
    \put(165,110){\makebox(0,0){$\bg{\gg}$}}

\end{picture}
\end{figure}
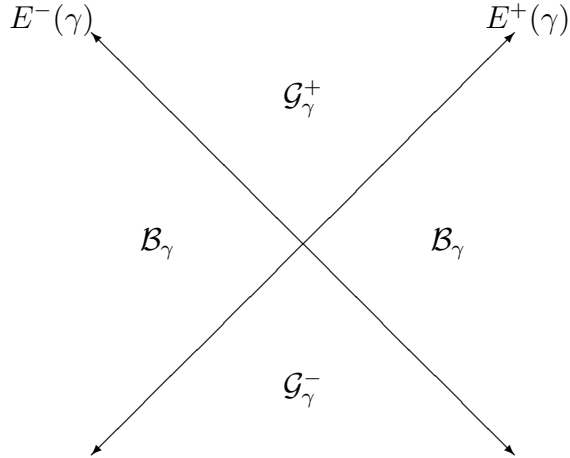

The closure of the two remaining components of $\E$:
\begin{align*} \ggp{\gg}&=\{ p_-\xm{\gg}+p_+\xp{\gg}+p_o\xo{\gg}
~|~ p_-,p_+\geq 0\}\\ \ggm{\gg}&=\{
p_-\xm{\gg}+p_+\xp{\gg}+p_o\xo{\gg} ~|~ p_-,p_+\leq 0\} ,
\end{align*} form $\cF(\gg)$. See Figure~\ref{fig:crosssection}.
The set $\ggp{\gg}$  is future complete.

Note that $\cF(\gg)=\cF(\gg^n)$, for every $n$, since
$\xpm{\gg^n}=\xpm{\gg}$. Thus: \begin{prop}\label{prop:fixedhyper}
Let $\gg$ be a hyperbolic transformation with fixed points. Then
\begin{equation*} \cF(\langle\gg\rangle)=\ggp{\gg}\cup
\ggm{\gg}\setminus C_{\gg} \end{equation*} and \begin{equation*}
\cT(\langle\gg\rangle)=\bg{\gg}. \end{equation*} \end{prop} The
interior of $\ggp{\gg}$ is a maximal connected open
$\gg$-invariant CTC-free region of $\E$ on which
$\langle\gg\rangle$ acts properly discontinuously.  The quotient
$\ggp{\gg}/\langle\gg\rangle$ is called {\em Misner space} and is
a  Lorentz spacetime which is future complete and diffeomorphic to
$\R^2\times S^1$.

The interior of $\ggm{\gg}$ yields an analogous, past-complete
spacetime.

\subsection{Hyperbolic transformations without fixed points}
Define $\bg{\gg}$, $\ggp{\gg}$, $\ggm{\gg}$ as above.  Solutions
to \eqref{eq:hyperbolicspacelike} now consist of two components
which are strictly contained in $\bg{\gg}$.  They are bounded by
hyperbolic sheets forming $\cL(\gg)$, which are asymptotic to the
planes $E^{\pm}(\gg)$.

\begin{thm}
Let $\gg\in\Isom(\E)$ be a hyperbolic isometry without fixed
points; then \begin{equation*} \cT(\langle\gg\rangle) = \bg{\gg},
\end{equation*}
and the CTC region of $\E/\langle\gg\rangle$ is $
\bg{\gg}/\langle\gg\rangle$. Furthermore, \begin{equation*}
\cF(\langle\gg\rangle)= \ggp{\gg}\cup\ggm{\gg}. \end{equation*}
The regions $\ggp{\gg}$ and $\ggm{\gg}$ are future complete and
past complete, respectively. All closed curves in
$\ggp{\gg}/\langle\gg\rangle$ and $\ggm{\gg}/\langle\gg\rangle$
are spacelike. \end{thm} \proof  Note that the linear part of
$\gg^n$ is $g^n$ so $\gl(g^n)= (\gl(g))^n$, but $\ga(\gg^n)=
n\ga(\gg)$. Consider the right hand side of
\eqref{eq:hyperbolicspacelike}. We find: \begin{equation*}
\lim_{n\rightarrow\infty}\frac{-\ga(\gg^n)^2 }{2
(1-\gl(g^n))(\gl(g^n)^{-1} -1)}=0. \end{equation*} Therefore, the
hyperbolic sheets $\cL(\gg^n)$ approach $E^{\pm}(\gg)$, since
$E^{\pm}(\gg^n)=E^{\pm}(\gg)$. \cqfd

In contrast to hyperbolic transformations with fixed points, the
maximal connected open subset of $\E$ on which $\langle\gg\rangle$
acts properly discontinuously is all of $\E$. The quotient
$\E/\langle\gg\rangle$ is a Lorentz spacetime which is future
complete, past complete, and diffeomorphic to $\R^2\times S^1$,
which Grant \cite{Grant} identified as representing the complement
of two straight moving cosmic strings which do not intersect.

Now that we understand the CTC regions for a group generated by a
single hyperbolic transformation, we can start to look at more
complicated groups. Here is an interesting example, where the CTC
region consists of all of $\E$.
\begin{ex}[$\cF(\gG) = \emptyset$]

Suppose that  $\gg_1$ and $\gg_2$ are hyperbolic transformations
such that \begin{itemize} \item the invariant lines $C_{\gg_1} ,
C_{\gg_2}$ are distinct, \item the  linear parts admit the same
eigenvectors, and \item the group $\gG = \langle \gg_1 , \gg_2
\rangle$ acts freely and properly discontinuously on $\E$.
\end{itemize}

Note that there are such groups. One such example is the group
generated by the transformations:
\begin{equation*}
    \begin{array}{rcl} \gg_1 (x) & = & \begin{bmatrix} 1 & 0 & 0 \\
                            0 & \frac{3}{2} & -\frac{\sqrt{5}}{2} \\
                            0 & -\frac{\sqrt{5}}{2}&  \frac{3}{2}
            \end{bmatrix} x +
            \begin{bmatrix} 1 \\ 0\\ 0 \end{bmatrix} \mbox{ and} \\ \gg_2 (x) & = & \begin{bmatrix} 1 & 0 & 0 \\
                            0 & \frac{7}{2} & -\frac{3\sqrt{5}}{2} \\
                            0 & -\frac{3\sqrt{5}}{2}&  \frac{7}{2}
            \end{bmatrix} x +
            \begin{bmatrix} 2 \\ \frac{1}{\sqrt{5}} \\ 1  \end{bmatrix}
    \end{array}
\end{equation*}
These transformations admit distinct invariant lines which are
both parallel to the vector$\begin{bmatrix} 1 \\ 0 \\0
\end{bmatrix}$. (In fact, $\E/\langle \gg_1 ,\gg_2 \rangle $ is a
$3$--torus. There is a normal subgroup $\cong \mathbb{Z}^2$ of
pure translations and the linear parts preserve this
$\mathbb{Z}^2$--lattice.)

To study $\cT(\gG )$, it is enough to look at a single cross
section parallel to
\begin{equation*}
    \langle \vx^-(g_1),
\vx^+(g_1) \rangle =\langle \vx^-(g_2), \vx^+(g_2) \rangle.
\end{equation*}
Conjugate invariant lines are obtained as follows:
\begin{equation*} C_{\gg_i^{} \gg_j^{} \gg_i^{-1}}= \gg_i C_{\gg_j}.
\end{equation*} Suppose that $q\in\ggp{\gg_1}\cup \ggp{\gg_2}$.
The  point representing the invariant line $C_{\gg_1^n \gg_2^{}
\gg_1^{-n}}$ approaches the asymptote representing either
$E^{+}(\gg_1)$ or $E^{-}(\gg_1)$ as seen in
Figures~\ref{fig:bad1},\ref{fig:bad2} and \ref{fig:bad5}. So in
our example, every point $q$ lies in $\bg{\gg_1^n \gg_2^{}
\gg_1^{-n}}$ for a large enough $n$. The other cases are treated
in the same manner, showing that every $q$ lies in $\bg{\gg}$, for
some $\gg\in\gG$.
\begin{figure}[ht]
\centerline{ \epsfclipon \epsfxsize=6cm \epsfbox{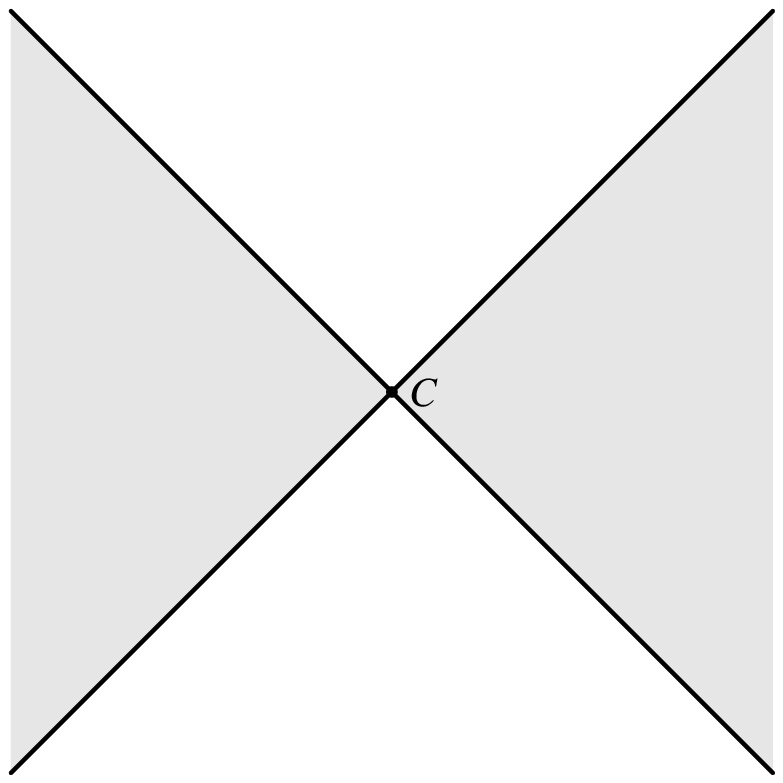}}
    \caption{A cross-section of $\bg{\gg_1^{}}$; $C=C_{\gg_1^{}}$.
    }
    \label{fig:bad1}
\end{figure}

\begin{figure}[ht]
    \centerline{ \epsfclipon \epsfxsize=6cm \epsfbox{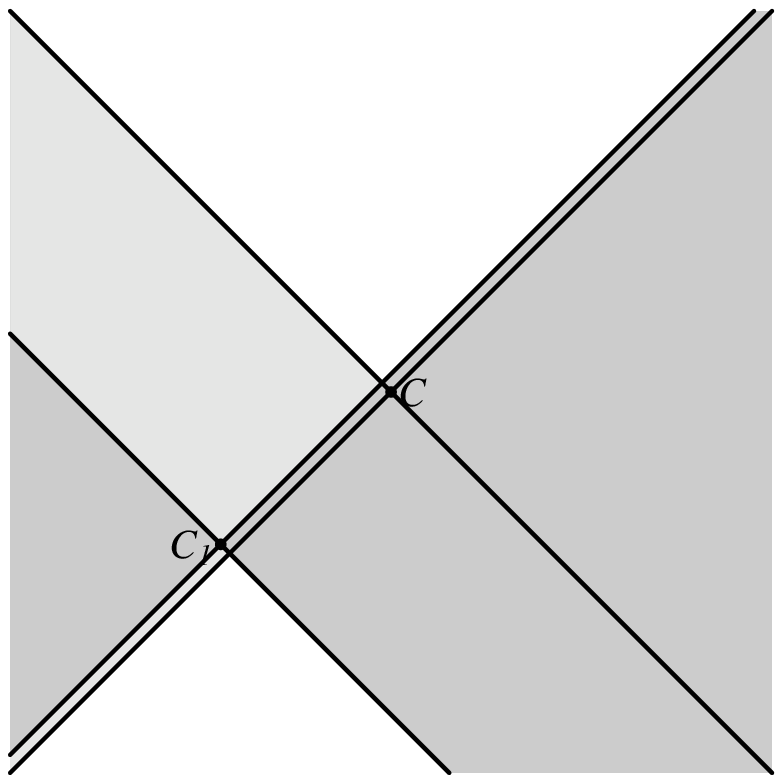}}
    \caption{A cross-section of 
    $\bg{\gg_1^{}}$; $C_1=C_{\gg_2^{}}$
    }
    \label{fig:bad2}
\end{figure}

\begin{figure}[ht]
    \centerline{ \epsfclipon \epsfxsize=6cm \epsfbox{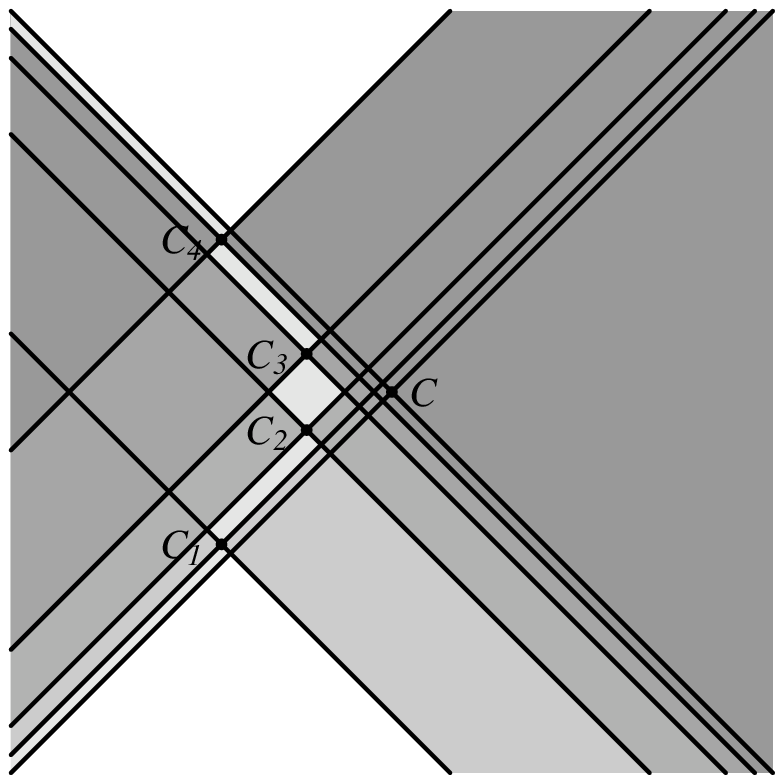}}
    \caption{A cross-section of $\bg{\gg_{1}^i\gg_2^{}\gg_{1}^{-i}}$,
    $0\leq i\leq 3$, and $\bg{\gg_1^{}}$;
    $C_i=C_{\gg_{1}^{i-1}\gg_2^{}\gg_{1}^{-i+1}}$
    }
    \label{fig:bad5}
\end{figure}

Thus, the CTC region  of $\E/\Gamma$ is the entire space.
\end{ex}

\section{Parabolic transformations}
In the linear group $G$, all parabolic transformations are
conjugate to each other. The conjugacy classes of affine parabolic
transformations in $\Isom(\E)$ are described in
\cite{CharetteDrumm}. In particular, the invariant $\ga$, which is
defined for hyperbolic transformations only, can be generalized to
parabolic transformations.

\subsection{Parabolic transformations with fixed points}
Let $\rho\in\Isom(\E)$ be a parabolic isometry with a fixed point.
Equivalently, $\rho$ admits a pointwise fixed lightlike line,
which is parallel to the unique (fixed) eigendirection of $\rho$'s
linear part.

For simplicity, identify $\rho$ with $\rho$'s linear part in $G$
-- this may be achieved by choosing an origin $\origin$ in the
fixed point set of $\rho$.  (We will omit $\origin$ in the rest of
this discussion.)  Thus $\rho$ admits a fixed eigendirection of
lightlike vectors.  Choose one such vector which is future
pointing and call it $\xo{\rho}$.

Let $\xo{\rho}^\perp$ denote its Lorentz-orthogonal plane in
$\R^{2,1}$, that is, the set of vectors whose Lorentz inner
product with $\xo{\rho}$ vanishes.  Note that $\xo{\rho}^\perp$ is
tangent to the light cone and contains $\R\xo{\rho}$, the line
spanned by $\xo{\rho}$.  Every vector in $\xo{\rho}^\perp$ is
either lightlike (if and only if it is parallel to $\xo{\rho}$) or
spacelike; in particular, $\xo{\rho}^\perp$ contains no timelike
vectors.

Similarly to the hyperbolic case, let $\vx^1$, $\vx^2$ be a pair
of vectors such that:
\begin{itemize}
    \item $\vx^1$ lies in $\xo{\rho}^\perp$ and is unit-spacelike;
    \item $\vx^2$ is lightlike, future pointing and Lorentz-orthogonal
    to $\vx^1$, but not parallel to $\xo{\rho}$;
    \item $\{\xo{\rho},\vx^1,\vx^2\}$ is a positively oriented basis.
\end{itemize}
Then relative to this basis, every power of $\rho$ can be written as
an upper-triangular matrix with 1's on the diagonal:
\begin{equation*} \rho^n = \begin{bmatrix}
          1 & a_n & b_n \\
          0 & 1   & c_n \\
          0 & 0   & 1
         \end{bmatrix}.
\end{equation*}
In particular, for every $p\in\E$ each displacement vector
$\rho^n(p)-p$ lies in $\xo{\rho}^\perp$. (More generally, if
$\gg\in\Isom(\E)$ admits fixed points, then the displacement
vectors lie in $\xo{\rho}^\perp$. We will use this fact again in
the elliptic case.) Thus the displacement vector is never
timelike.  The displacement vector is lightlike if and only if
$p\in\xo{\rho}^\perp$  and, of course, vanishes if and only if $p$
is fixed by $\rho$.

Thus:
\begin{align*}
    \cS(\rho)&=\cS(\rho^n)=\E -\xo{\rho}^\perp\\
    \cL(\rho)&=\cL(\rho^n)=\xo{\rho}^\perp .
\end{align*} We define
\begin{equation*}
    \E_{\rho} = \E - \R\xo{ \rho } ,
\end{equation*}
and we
have shown the following:
\begin{prop} Let $\rho$ be a parabolic
isometry with fixed points.Then:
\begin{equation*}
\cF(\langle\rho\rangle)=\E_\rho \mbox{ and } \cT(\langle \gg
\rangle )= \emptyset .
\end{equation*}
That is, the CTC region of
$\E_\rho /\langle\rho\rangle$ is empty. The spacelike region
$\cS(\rho)$ divides into two components separated by the lightlike
region $\cL(\rho)=\xo{\rho}^\perp - \R\xo{\rho}$.
\end{prop}

\subsection{Parabolic transformations with no fixed points}
The situation changes dramatically for parabolic transformations
without fixed points.

Recall that all linear parabolics are conjugate to each other.
Choose $\rho$ such that: \begin{equation*}
\{\xo{\rho},\vx^1,\vx^2\} = \left\{ \begin{bmatrix} 0 \\
\frac{-1}{\sqrt{2}} \\ \frac{1}{\sqrt{2}} \end{bmatrix},
\begin{bmatrix} 1 \\ 0 \\ 0 \end{bmatrix}, \begin{bmatrix} 0 \\
\frac{1}{\sqrt{2}} \\ \frac{1}{\sqrt{2}} \end{bmatrix} \right\}
.\end{equation*} Choose an origin $\origin$ and denote the
coordinates of a point $p$ in the $\{\xo{\rho},\vx^1,\vx^2\}$
basis as follows: \begin{equation*} p=\begin{bmatrix}\tx \\ \ty \\
\tz \end{bmatrix}. \end{equation*} We may further conjugate the
parabolic isometry so that its translational part is parallel to
$\vx^2$; then the transformation may be written:
\begin{equation}\label{eq:parabolicewithout}
\rho_{\gt}(p) = \begin{bmatrix} 1 & \sqrt{2} & 1 \\
            0 & 1 & \sqrt{2} \\
            0 & 0 & 1
          \end{bmatrix}
          \begin{bmatrix} \tx \\ \ty \\ \tz  \end{bmatrix}+
          \begin{bmatrix} 0 \\ 0 \\ \gt \end{bmatrix}.
\end{equation}
          The inner product in this basis is given by
\begin{equation*}
    \B ( p ,q ) = -\tx\tw + \ty \tv - \tz \tu ,
\end{equation*}
where $ \tu ,\tv , \tw$ are the coordinates of $q$ in the new
basis. A point $p$ is  in $\cT(\rho_{\gt})$ if:
\begin{equation*}
    \tz^2 - \gt (\sqrt{2}\ty +\tz) < 0 .
\end{equation*}
More explicitly,
\begin{equation}\label{eq:parabolicsheet}
    \left\{ \begin{array}{lcr}
    \ty > \frac{\tz^2-\gt\tz}{\sqrt{2}\gt} & \mbox{ if } & \gt > 0 \\
    \ty < \frac{\tz^2-\gt\tz}{\sqrt{2}\gt} & \mbox{ if } & \gt < 0
    \end{array}
    \right.
\end{equation}
We see that $\cL(\rho_{\gt})$ is a parabolic sheet bounding
$\cT(\rho_{\gt})$. Let us examine the case where $\gt>0$. Suppose
that $p\in \cL(\rho_{\gt}) \cup \cT (\rho_{\gt})$. The set of
future pointing, non-spacelike vectors is convex. Therefore,
\begin{equation*} \begin{array}{rcl} \rho^{n+1}_{\gt}(p) -p & = &
\left( \rho^{n+1}_{\gt}(p)-
            \rho_{\gt}(p)\right) +
            \left( \rho_{\gt}(p)- p \right)\\
 & = & g\left( \rho_{\gt}^n(p)- p \right)+
            \left( \rho_{\gt}(p)- p \right)
\end{array}
\end{equation*}
is a future pointing timelike vector, since each term is either a
future pointing lightlike vector or a future pointing timelike
vector and all vectors are not parallel.

The same statement holds when $\gt <0$, substituting the term
``past'' for ``future''.

\begin{lemma}\label{timelikeinclusion}
Let $\gt\neq 0$; then
$\cT(\rho_{\gt}^n)\subset\cT(\rho_{\gt}^{n+1})$ for all $n\geq 1$.
\end{lemma}

We will now show by direct calculation that the parabolic sheet
$\cL(\rho_{\gt}^n)$ is described by the equation: \begin{equation}
\ty = \frac{\tz^2-\gt\tz}{\sqrt{2}\gt}-\phi(n)\gt , \end{equation}
where $\phi(n)$ increases as fast as $n^2$. Notice that this sheet
is just a translate of $\cL(\rho_{\gt})$ in the $\vx^1$ direction
by the amount $-\phi(n)\gt$.

We can show the following by induction:
\begin{equation}
\rho_\gt^n(p)= \left[\begin{array}{c}
n\sqrt{2}\ty+n^2\tz+\gt\sum_{i=1}^{n-1}i^2 \\ n\sqrt{2}\tz
+\sqrt{2}\gt\sum_{i=1}^{n-1}i\\ n\gt
\end{array}\right] .
\end{equation}

We get an equation describing $\cL(\rho_{\gt}^n)$ by solving
\begin{equation*}
\B\left(\rho_\gt^n(p)-p,\rho_\gt^n(p)-p\right)=0,
\end{equation*}
and we obtain
\begin{equation*}
n^2\sqrt{2}\gt\ty=n^2\tz^2+n\gt\tz\left(-n^2+2\sum_{i=1}^{n-1}i\right)
-\gt^2\left(n\sum_{i=1}^{n-1}i^2-\left(
\sum_{i=1}^{n-1}i\right)^{2}\right).
\end{equation*}
We note that $-n^2+2\sum_{i=1}^{n-1}i=-n$ and that
$6\sum_{i=1}^{n-1}i^2 = n(n-1)(2n-1)$, so the following holds:
\begin{equation*} \ty=\frac{\tz^2-\gt\tz}{\sqrt{2}\tau}
-\frac{\gt(n^2-1)}{12\sqrt{2}}.
\end{equation*}

\begin{thm}\label{thm:ParabolicNoFixed}
Let $\rho\in\Isom(\E)$ be a parabolic isometry without fixed
points.  Then: \begin{equation*} \cT(\langle\gg\rangle) = \E
\end{equation*} and the CTC region of $\E/\langle\gg\rangle$ is
the entire space. \end{thm}

\begin{ex}
Suppose $\{g_n\}_{n\geq 0}\subset G$ is a sequence of hyperbolic
isometries with a common fixed point $\origin$, converging to a
parabolic element in $G$.  This happens, for instance, by letting
$g_n=h^ngh^{-n}$, where $g$ is hyperbolic and $h$ is an arbitrary
element of $G$.  Then: \begin{equation*} \cT (g_n) \rightarrow \cT
(\rho)\mbox{ as } g_n \rightarrow \rho . \end{equation*} Indeed,
as $g_n\rightarrow\rho$, the stable and unstable planes of $g_n$
both approach the plane $\xo{\rho}^\perp$. (Recall that
$\xo{\rho}^\perp$ is tangent to the light cone at the origin and
contains $\R\xo{\rho}$.)  Thus $\cT(\gg_n)$ approaches the empty
set.

Now consider a sequence of hyperbolic transformations
$\{\gg_n\}_{n\geq 0}$ approaching a parabolic transformation
$\rho$, without fixed points. The regions $\cT(\gg_n)$ still
approach the empty set. However, $\cT(\rho)=\E$. \end{ex}

\section{Elliptic transformations}
In $G$, all elliptic elements are conjugate to an element of the
form \begin{equation}\label{eq:elliptic}
    \psi_{\theta} = \begin{bmatrix} \cos\theta & \sin\theta & 0 \\
                                    -\sin\theta & \cos\theta & 0 \\
        0 & & 1
    \end{bmatrix}.
\end{equation}
The signed Lorentzian length invariant $\ga$ cannot be generalized
in any coherent manner to elliptic elements. We will first
consider elliptic transformations with fixed points and then
without fixed points.

\subsection{Elliptic transformations with fixed points}
After choosing an origin $\origin$ in the fixed point set, we let:
\begin{equation*} p = \begin{bmatrix} x \\ y \\ z \end{bmatrix}
\end{equation*} so that the transformation can be written as
$\psi_{\theta}$. A fixed eigenvector for this transformation is
\begin{equation*} \xo{ \psi_{ \theta } }=\begin{bmatrix} 0 \\ 0 \\
1  \end{bmatrix} \end{equation*} and the fixed point set for
$\psi_{\theta}$ consists of the line $\R\xo{\psi_{\theta}}$.

If $\theta$ is a rational multiple of $2\pi$, then the group
$\langle \psi_{\theta} \rangle$ acts freely and properly
discontinuously on the complement of its fixed point set:
\begin{equation*} \E_{\psi_{\theta}} = \E-\R\xo{\psi_{\theta} }.
\end{equation*}

Recall that since $\xo{\psi_{\theta}}$ is fixed by
$\psi_{\theta}$, every displacement vector lies in the
Lorentz-orthogonal plane of $\xo{\psi_{\theta}}$.  This is the
$xy$-plane, which is spacelike. Thus: \begin{prop} Suppose $\psi$
is an elliptic isometry with fixed points, and that $\psi$ is a
rotation of a rational multiple of $2\pi$ about its line of fixed
points.  Then: \begin{equation*}
\cF(\langle\psi_{\theta}\rangle)=\E_{\psi_{\theta}} \mbox{ and }
\cT(\langle \psi_{\theta} \rangle )= \emptyset . \end{equation*}
That is, the CTC region of $\E_{\psi_{\theta}} /\langle
\psi_{\theta} \rangle$ is empty. \end{prop}

Spaces of the type $\E_{\psi_\theta}/\langle \psi_\theta \rangle$
can be identified with special cases of spacetimes that represent
(spinless) particles in 2+1-dimensional gravity.  Such spaces have
been described by Deser, Jackiw, and 't Hooft~\cite{DJTH}.  (See
also \cite{Welling} and the references cited there.) The direct
metric product of such a space with a spacelike line is known as a
cosmic string.

\subsection{Elliptic transformations without fixed points}
As in the  parabolic case, the situation changes dramatically for
elliptic transformations without fixed points. Every elliptic
transformation is conjugate to a transformation
\begin{equation}\label{eq:ellipticnofixedpoints}
\psi_{\theta,t}(p) = \begin{bmatrix} \cos\theta & \sin\theta & 0 \\
                                    -\sin\theta & \cos\theta & 0 \\
        0 & & 1
    \end{bmatrix}
          \begin{bmatrix} x \\ y \\ z \end{bmatrix}+
            \begin{bmatrix} 0 \\  0\\ t \end{bmatrix},
\end{equation}
which we will write as $\psi$. The group $\langle \psi\rangle$ acts
properly discontinuously on all of $\E$.

\begin{thm}\label{thm:EllipticNoFixed}
Let $\psi\in\Isom(\E)$ be an elliptic isometry without fixed
points.  Then: \begin{equation*} \cT(\langle\psi\rangle) = \E
\end{equation*}
and the CTC region of $\E/\langle\psi\rangle$ is the entire space.
\end{thm} \proof
    Note that
    \begin{equation*}
 \psi^{k} = \begin{bmatrix} \cos k\theta & \sin k\theta & 0 \\
                                    -\sin k\theta & \cos k\theta & 0 \\
        0 & & 1
    \end{bmatrix}
          \begin{bmatrix} x \\ y \\ z \end{bmatrix}+
            \begin{bmatrix} 0 \\  0\\ kt \end{bmatrix},
    \end{equation*}
so that the length (equivalent to the underlying topology) of the
projection of $\psi(p)-p$ onto the $(x,y)$--plane is bounded by
$2\| x^{2} + y^{2} \|$.  However, the projection of $\psi(p)-p$
onto the $z$--axis is unbounded as $k\rightarrow\infty$. Thus, for
a sufficiently large power $k$, depending on the distance from $p$
to the $z$--axis, the vector $\psi^{k}(p)-p$ is timelike. \cqfd

\section{Future directions (pun intended)}
In a future note, we will look for closed timelike curves in
$X/\Gamma$ where $X\subset\E$ and $\Gamma$ is more complicated.
Some examples are given below.

For the case $X=\E$ and free $\Gamma$, we call $\E/\Gamma$ a {\em
Margulis spacetime}.  Theorem~\ref{thm:ParabolicNoFixed} has an
interesting and immediate consequence for Margulis spacetimes. In
\cite{Drumm}, Margulis spacetimes with noncyclic free fundamental
groups containing parabolic transformations were constructed.
Parabolic transformations  were shown to be very much like
hyperbolic transformations for questions concerning proper actions
of a group on  $\E$ in \cite{Drumm} (and \cite{CharetteDrumm} for
that matter). But we see here  that for questions concerning
closed timelike curves, the difference between hyperbolic and
parabolic transformations is tremendous.

We will also be keenly interested in surface groups, groups
isomorphic to the fundmental group of a closed surface. As shown
in \cite{Mess} and \cite{GM} these groups do not act properly
discontinuously on $\E$. However, Mess showed \cite{Mess} that
surface groups can act properly discontinuously on some subset
$X\subset\E$.

\makeatletter

\renewcommand{\@biblabel}[1]{\hfill#1.}\makeatother

\end{document}